\documentclass[10pt,a4paper]{amsart}

\usepackage[latin1]{inputenc}
\usepackage[T1]{fontenc}
\usepackage{eucal,url,amssymb,bm,stmaryrd,booktabs,enumerate}
\usepackage[pagebackref]{hyperref}
\usepackage{accents}

\makeatletter
\def\widebar{\accentset{{\cc@style\underline{\mskip10mu}}}}
\def\Widebar{\accentset{{\cc@style\underline{\mskip8mu}}}}
\makeatother

\theoremstyle{plain}
\newtheorem{theorem}{Theorem}[section]
\newtheorem{mainthm}{Theorem}

\newtheorem{lemma}[theorem]{Lemma}
\newtheorem{corollary}[theorem]{Corollary}

\theoremstyle{definition}

\theoremstyle{remark}
\newtheorem{remark}{Remark}[section]

\newtheorem*{acknowledgements}{Acknowledgements}

\numberwithin{equation}{section}

\newcommand{\LC}{\nabla}
\newcommand{\lc}{\nabla}
\newcommand{\RC}{\textup{R}}

\newcommand{\Ric}{\textup{ric}}
\newcommand{\ric}{\Ric}

\newcommand{\abs}[1]{\left\lvert #1\right\rvert}

\begin{document}
\title[Riemannian products which are conformally
  Einstein]{Riemannian products which are conformally equivalent to
  Einstein metrics}

\date{23 May 2008}

\author{Richard Cleyton}
\address[Cleyton]{Humboldt-Universität zu Berlin\\
Institut für Mathematik\\
Unter den Linden 6\\
D-10099 Berlin\\
Germany}

\email{cleyton@mathematik.hu-berlin.de}

\begin{abstract}
  Necessary and sufficient conditions for a Riemannian product to be
  conformally equivalent to an Einstein manifold are given.  Such spaces
  which are complete are characterized.
\end{abstract}

\maketitle

\section{Introduction}

The study of Riemannian metrics that are conformally related to
Einstein metrics go back to Brinkmann~\cite{MR1512193} and before and
has recently received attention by Gover \& Nurowski~\cite{MR2171895}
and Listing~\cite{MR2171892}.  It is subtly related to the subject of
twistor spinors, see e.g.~\cite{twistor-killing} and to the famous
Yamabe problem of whether the conformal class of a metric has
representative with constant scalar curvature.  Recently, Moroianu and
Ornea considered the following.  Take the product metric $g=g'+dt^2$
on $M=M'\times\mathbb R$ with $M'$ compact but not the round sphere.
Then $g$ is conformally equivalent to an Einstein metric $\bar g
=\phi^{-2}g$ if and only if the conformal factor $\phi$ is essentially
$\cosh(t)$.  Here ``essentially'' means that we identify homothetic
metrics on $M$ and $M'$ and ignore translations of the parameter $t$.
In particular $\phi$ does not depend on the coordinate in $M'$.  This
result was generalized by J.~M.~Ruiz even more recently to the case of
$M=M'\times \mathbb R^q$ equipped with the product metric which is the
standard metric on the second factor.  He shows amongst other things
that for $q>1$ no conformal factor will change such a metric to a
scalar positive Einstein metric~\cite{ruiz-2008}.  These results may
be generalized even further.  The reader may also wish to
consult~\cite{MR608388} and ~\cite{MR667635} for some related results.

\begin{mainthm}\label{thm:1}
  Suppose $(M^n,g)$ is a Riemannian product. Suppose $\phi$ is a
  positive function on $M$ such that $\bar g:=\phi^{-2}g$ is Einstein
  with Einstein constant $\bar\lambda$.
  In this case then either
  \begin{enumerate}[{\rm (I)}]
  \item $\bar g= g_B + f^2 g_F$ is a warped product metric on $M =
    B^q\times F^p$ and $g=f^{-2}g_B+g_F$ where $f\colon B\to
    \mathbb R$ is a positive function such that
    $\ric(g_F)=\lambda_Fg_F$,
    $\Delta_{g_B}f-(p-1)\abs{df}^2+\lambda_F=\bar\lambda f^2$,
    $f^2\ric(g_B) = pf \lc^{g_B}df +\bar\lambda g_B$.\label{item:1}
  \item $M^n=M^{n_1}_1\times M_2^{n_2}$ and $g=g_1+g_2$, where both
    $g_1$ and $g_2$ are Einstein metrics with Einstein
      constants $\lambda_1$ and $\lambda_2$.  Furthermore, there are
    non-constant functions $\phi_1$ on $M_1$ and $\phi_2$ on $M_2$
    such that
    \begin{enumerate}[{\rm (a)}]
    \item $\phi(x,y)=\phi_1(x)+\phi_2(y)$ for all $(x,y)$ in $M$
    \item constants $\bar a$ and $\bar b$ exists to the effect that
      $\lc^{g_1}d\phi_1 = -(\bar a\phi_1 + \bar b)g_1$,
      $\lc^{g_2}d\phi_2 = (\bar a\phi_1 - \bar b)g_2$, $\lambda_1 =
      (n_1 - 1)\bar a$ and $\lambda_2 = -(n_2 - 1)\bar a$.
    \end{enumerate}
    \label{item:2}
  \end{enumerate}
  Conversely, if $(M,g)$ and $\phi$ fulfill the conditions of either
  \eqref{item:1} or~\eqref{item:2} then $\phi^{-2}g$ is Einstein. When
  the conditions of \eqref{item:2} holds, $\bar c_1 :=
  \abs{d\phi_1}^2+\bar a\phi_1^2+2\bar b\phi_1$ and $\bar c_2 :=
  \abs{d\phi_2}^2 - \bar a\phi_2^2 + 2\bar b\phi_2$ are constants such
  that $c_1+c_2 = \bar \rho$, the renormalized scalar curvature of
  $\bar g$. 
\end{mainthm}
The conditions of \eqref{item:1} are precisely those giving an
Einstein warped product metric, see \cite{Besse:Einstein}.  On the
other hand the conditions of \eqref{item:2} force both $g_1$ and $g_2$
to be warped products of a very restricted form, at least locally.  An
argument in the proof of Lemma 5.2 of~\cite{math.DG/0610599} shows
that $g$ is also, locally, a warped product.  Note that, if $B$ is
compact the metrics $\bar g$ of \ref{item:1} have non-negative scalar
curvature.  The metrics under \ref{item:2} come with both positive as
well as negative scalar curvature.

In the situation of a complete metric $g$ we can be more concrete.
\begin{mainthm}\label{thm:2}
  Suppose $(M,g)$ is a Riemannian product and $g$ is complete.  If $g$
  is conformally equivalent to an Einstein metric $\bar g$ which is
  not a warped product, then $M$ is $\mathbb R^n = \mathbb
  R^{n_1}\times\mathbb R^{n_2}$ with the standard metric $g=g_1+g_2$.
  Using polar coordinates such that $g_1=ds^2+s^2g_{S^{n_1-1}}$ and
  $g_2=dt^2+t^2g_{S^{n_2-1}}$ the conformal factor is
  $\phi(s,t)=\tfrac12(s^2+t^2+R^2)$ for some positive constant $R$.
\end{mainthm}
\emph{Note:} Setting $r\cos\theta=s,~r\sin\theta=t$ the metric $g$
becomes $g = dr^2 + r^2g^*$ where $g^* := d\theta^2 + \cos^2\theta
g_{S^{n_1-1}} + \sin^2\theta g_{S^{n_2-1}}$.  We note that, where
defined, $g^*$ agrees with the standard metric of $S^{n-1}$.  In this
way we recognize $g$ as the standard metric on $\mathbb R^n$.  Under
the transformation $(s,t)\mapsto(r,\theta)$ the conformal factor
becomes $\tfrac12(r^2 + R^2)$.  Suitably restricted, the original
metric $g$ is a warped product metrics on the half-line times the
$(n-1)$-sphere.  Under the map $u=2\tan^{-1}(r/R)/R$, the conformally
changed metric $\bar g$ goes to $\bar g = du^2+\sin^2(u)g_{S^{n-1}}$,
the standard metric of the $n$-sphere.  In particular, $\bar g$ has
positive scalar curvature.  However, $\bar g$ is not complete as the
conformal transformation $u(r)$ isn't onto.  More generally,
part~\ref{item:2} provide instances of conformal changes which change
the sign of the scalar curvature (see the examples
of~\cite{MR979791}).  When, under~\ref{item:1}, $B$ is compact it is
known that the conformally changed metric is either scalar positive or
$f$ is constant, see~\cite{MR1974657}.

\subsection*{Some comments and consequences}

Different variations along the lines of warped products have been made
in literature.  One such is the \emph{doubly warped product metric}
$\bar g=b^2 g_B+ f^2 g_F$ with $b\colon F\to\mathbb R^+$, $f\colon B
\to\mathbb R^+$.  We note that $g := (fb)^{-2} \bar g$ is a Riemannian
product.  Thus $g$ is mixed Ricci-flat if and only if functions
$\phi_1,\phi_2$ exist such that $\phi_1(x) + \phi_2(y) = 1/f(x)b(y)$
whence $X(f)Y(b)=0$ for all vector field $X$ on $B$,~ $Y$ on $F$.  In
particular a doubly warped product is Einstein only if it is a warped
product.  A \emph{twisted product} is a metric $\bar g = g_B + f^2
g_F$ where $f$ is a positive function defined on $B\times F$.  When
the dimension of $F$ is greater than $1$ such metrics are also known
to be mixed Ricci-flat only in the case of $f$ being a product:
$f(x,y)=f_1(x)f_2(y)$~\cite{MR1865565}.  In particular $\bar g$ is
conformal to the Riemannian product $f_1^{-2}g_B+f_2^2g_F$ and whence
twisted products are Einstein only if they are warped.  \emph{Doubly
  twisted} metrics are $b^2g_B+f^2g_F$ with $b,g$ being positive
functions on the product $B\times F$.  Doubly twisted metrics include
metrics conformally equivalent to Riemannian products as a special
case.  As a corollary of our discussion above we have

\begin{corollary}
  There exists Einstein metrics which are doubly twisted but not warped.
\end{corollary}

Comparing the claims we make in Theorems~\ref{thm:1} and \ref{thm:2}
with the investigations of Moroianu and Ornea another issue
arises.  Is it possible, in general, to characterize the geometry of
$(M_1^p,g_1)$ and $(M_2^{n-p},g_2)$ such that the product $M^n$ is
either conformally equivalent to a (strict) nearly K\"ahler metric in
dimension $6$ or a $7$ dimensional weak holonomy $G_2$ manifold?
These two types of geometries are Einstein with positive scalar
curvature and, both before and after conformal change, have a nice
characterization in terms of Fernand\'ez-Gray-Hervella
classes~\cite{math.DG/0607487}.  They are also spaces with non-trivial
Killing spinors and so are conformal precisely to spaces admitting
twistor spinors.  In the case $p=1$ the only possible way to obtain a
Riemannian product conformal to a nearly K\"ahler $6$-fold is by
equipping $M_2$ with a Sasaki-Einstein metric and taking $M_1$ to be
an interval, see~\cite{math.DG/0602160, math.DG/0610599}.

\section{Preliminaries}

Let $(M^n,g)$ denote a Riemannian manifold of dimension $n\geqslant
3$.  Write $\LC$ for the Levi-Civita connection, $\RC^g$ for the
Riemannian curvature, $\Ric(g)$ for the Ricci curvature and $s=s(g)$
for the scalar curvature of the metric $g$.  The normalized scalar
curvature $\rho=s/n(n-1)$.  In what follows we shall be sloppy and
omit pull backs and inclusions so that on product manifolds $M_1\times
M_2$ a vector field $X$ on $M_1$ is also a vector field on the product
(since it defines one uniquely) and metrics $g_1$ and $g_2$ on the
factors define a unique metric, denoted $g_1+g_2$, on the product by
adding the pull backs of $g_1$ and $g_2$.  With these conventions the
Levi-Civita connection of $g$ is related to the Levi-Civita
connections on the factors by
\begin{gather}
  \LC^g_{X_i}Y_i=\LC^{g_i}_{X_i}Y_i,\qquad
  \LC^g_{X_i}Y_j=0,\quad\text{if $i\not=j$}. 
\end{gather}
for $X_1,Y_1\in\Gamma(TM_1)$ and $X_2,Y_2\in \Gamma(TM_2)$.  The
Riemannian curvature, the Ricci curvature and scalar curvature,
respectively, satisfy
\begin{gather}
  \label{eq:6} \RC^g = \RC^{g_1}+\RC^{g_2},\qquad \Ric(g) =
  \Ric(g_1)+\Ric(g_2),\qquad s(g)=s(g_1)+s(g_2).
\end{gather}

\subsection{Conformal change of Curvature}
Suppose $\phi\colon M\to\mathbb R^+$ is a smooth function.  Set $\bar
g:=\phi^{-2}g$.  This much is certainly well known:
\begin{equation}
  \label{eq:20}
  \LC^{\bar g}_XY = \LC_XY - \phi^{-1}(d\phi(X)Y+d\phi(Y)X) +
  \phi^{-1}g(X,Y)D^g\phi, 
\end{equation}
where $D^g\phi$ is the gradient of $\phi$, and
\begin{align}
  \label{eq:9}\Ric(\bar g) &= \Ric(g) + (n-2)\phi^{-1}\LC^g(d\phi) -
  (\phi^{-1}\Delta_g\phi + (n-1)\phi^{-2}\abs{d\phi}^2)g,\\
  \label{eq:21}
  \bar s &= \phi^2 s - 2(n-1)\phi\Delta_g\phi -n(n-1)\abs{d\phi}^2,
\end{align}
where $\Delta_g\phi$ denotes the Laplacian of $\phi$, see for
instance~\cite{Besse:Einstein}, Chapter~1.J.  A metric $g$ is
\emph{Einstein} if there is a constant $\lambda$, called \emph{the
  Einstein constant}, such that $\Ric(g) = \lambda g$.  An Einstein
metric with $\lambda=0$ is \emph{Ricci-flat}.  A metric is
\emph{conformally Einstein} if a positive function $\phi$ on $M$ and a
constant $\bar\lambda$ exists such that
\begin{equation}\label{eq:17}
  \Ric(g) = - (n-2)\phi^{-1}\LC^g(d\phi) + (\phi^{-2}\bar\lambda +
  \phi^{-1}\Delta_g\phi + (n-1)\phi^{-2}\abs{d\phi}^2)g. 
\end{equation}

It is well-known that the metric of a Riemannian product is Einstein
if and only if the factors are Einstein for the same Einstein
constant.  Applying a homothety to an Einstein metric gives another
Einstein metric.  If two Einstein metrics are conformally equivalent
then the conformal factor satisfies the second order differential
equation $\LC^gd\phi=-\tfrac1n \Delta_g\phi g$.

\section{Conformally Einstein Riemannian products}
Before dealing with Einstein metrics we shall prove a couple of
related but more general results.  First, by equation~\eqref{eq:6},
the Ricci curvature of a Riemannian product is a section of
$S^2(T^*M_1)\oplus S^2(T^*M_2)$, where $S^2$ is the second symmetric
power of the (pulled back) bundle.  In general, a metric $g$ on a
product manifold is \emph{mixed Ricci-flat}~\cite{MR1865565} if it has
the property that $\Ric(g)$ is a section of $S^2(T^*M_1)\oplus
S^2(T^*M_2)$. 
\begin{lemma}\label{lem:1}
  Suppose $(M,g)$ is a Riemannian product.  Let $\phi\colon
  M\to\mathbb R$ be a positive function.  Then $\bar g = \phi^{-2} g$
  is mixed Ricci-flat if and only if there exists functions
  $\phi_1\colon M_1\to \mathbb R$ and $\phi_2\colon M_2\to \mathbb R$
  such that $\phi(x,y) = \phi_1(x) + \phi_2(y)$ for all $(x,y)$ in
  $M$.
\end{lemma}
\begin{proof}
  The metric $\phi^{-2}(g_1 + g_2)$ is mixed Ricci-flat if and only if
  $0 = \Ric(\bar g)(X,Y) = \LC_{X}(d\phi)(Y) = X(Y(\phi))$ for all
  $X\in \Gamma(TM_1),~Y\in \Gamma(TM_2)$.  Clearly if $\phi(x,y) =
  \phi_1(x) + \phi_2(y)$ then this is satisfied.  Conversely, keep
  $x_0$ in $M_1$ fixed and consider $\psi(x,y) := \phi(x,y) -
  \phi(x_0,y)$.  Then $Y(\psi)(x,y) = Y(\phi)(x,y) - Y(\phi)(x_0,y) =
  0$ since $Y(\phi)$ is independent of $x$.  Therefore we may set
  $\phi_1(x)=\psi(x,y)$ and $\phi_2(y)=\phi(x_0,y)$.
\end{proof}
\begin{remark}
  The functions $\phi_i$ of Lemma~\ref{lem:1} will be called the
  summands of the conformal factor $\phi$.  They are unique only up to
  the transformation $(\phi_1,\phi_2)\to(\phi_1+c,\phi_2-c)$ where $c$
  is a constant.  Therefore, if say $\phi_2$ is constant, no
  generality is lost by setting $\phi_2=0$.
\end{remark}

\begin{lemma}\label{lem:2}
  Let $(M,g)$ be a Riemannian product with constant scalar
  curvature. Suppose $\bar g = \phi^{-2}g$ is mixed Ricci-flat and
  also has constant scalar curvature.  Let $\phi_1$ and $\phi_2$ be
  the summands as in Lemma \ref{lem:1} and let $\Delta_i$ be the
  Laplacians of the factors in $M$.  If both $\phi_1$ and $\phi_2$ are
  non-constant, then real numbers $a_1,b_1,c_1,a_2,b_2,c_2$ exist such
  that
  \begin{gather}
    \begin{cases}
      \Delta_1\phi_1=a_1\phi_1+b_1,\quad\Delta_2\phi_2=a_2\phi_2+b_2,\\
      c_1 = (a_2-a_1)\phi_1^2  - 2(b_1+b_2)\phi_1 - n\abs{d\phi_1}^2,\\
      c_2 = (a_1-a_2)\phi_2^2  - 2(b_1+b_2)\phi_2- n\abs{d\phi_2}^2,\\
      s = (n-1)(a_1+a_2),
    \end{cases}\label{eq:25}\\
    \bar s = (n-1)(c_1+c_2).\label{eq:26}
  \end{gather}
  Conversely, suppose that $(M_1,g_1)$ and $(M_2,g_2)$ are Riemannian
  manifolds with constant scalar curvatures $s_1$ and $s_2$, and that
  $\phi_1\colon M_1\to\mathbb R, ~\phi_2\colon M_2\to\mathbb R$ and
  $(a_1,a_2,b_1,b_2,c_1,c_2)\in\mathbb R^6$ solves \eqref{eq:25} with
  $s=s_1+s_2$.  Set $\phi:=\phi_1+\phi_2$, $M=M_1\times M_2$ and
  define $M'=M\setminus \{p\colon \phi(p)=0\}$.  Then $\bar g =
  \phi^{-2}(g_1+g_2)$ is a constant scalar curvature metric on $M'$
  with scalar curvature $\bar s$ given by~\eqref{eq:26}.
\end{lemma}
\begin{proof}
  Suppose that $\phi_1$ and $\phi_2$ both are non-constant.  Take $X$
  in $\Gamma(TM_1)$ and $Y$ in $\Gamma(TM_2)$ and
  differentiate~(\ref{eq:21}).  This gives
  \begin{equation}\label{eq:23}
    \frac{s}{n-1}X(\phi_1)Y(\phi_2)=X(\Delta_1\phi_1)Y(\phi_2)
    +
    X(\phi_1)Y(\Delta_2\phi_2).
  \end{equation}
  Write $Z:= D^{g_1}\phi_1,~W:=D^{g_2}\phi_2$ and set $U_1:=\{x\in
  M_1\colon Z_x\not=0\},~U_2:=\{y\in M_2\colon W_y\not=0\}$.  Then in
  the open subset $U_1\times U_2$ of $M$ we have
  \begin{equation*}
    \frac{s}{n-1} = \frac{Z(\Delta_{g_1}\phi_1)}{Z(\phi_1)} +
    \frac{W(\Delta_{g_2}\phi_2)}{W(\phi_2)} 
  \end{equation*}
  in $U_1\times U_2$.  Therefore there are constants $a_1$ and $a_2$
  with $(n-1)(a_1+a_2)= s$ such that
  $Z(\Delta_{g_1}\phi_1)=a_1Z(\phi_1)$ in $U_1$ and
  $W(\Delta_{g_2}\phi_2)=a_2W(\phi_2)$ in $U_2$.  But
  $Z(\Delta_{g_1}\phi_1)=a_1Z(\phi_1)$ holds trivially in the interior
  of $M_1\setminus U_1$ and so by continuity it holds throughout
  $M_1$.  Now insert $Z$ instead of $X$ in \eqref{eq:23}.  Then we
  have
  \begin{equation}
    (a_1 + a_2)Z(\phi_1)Y(\phi_2) = a_1Z(\phi_1)Y(\phi_2)
    + Z(\phi_1)Y(\Delta_{g_2}\phi_2). 
  \end{equation}
  and so $Y(\Delta_{g_2}\phi_2 - a_2\phi_2) = 0$ for all $Y$
  throughout $M_2$ and hence $\phi_2$ satisfies
  $\Delta_{g_2}\phi_2=a_2\phi_2+b_2$ for some constant $b_2$.  By
  symmetry, the analogous equation holds for $\phi_1$.  Inserting this
  in \eqref{eq:21} gives
  \begin{multline*}
    \frac{\bar s}{n-1}=\left( (a_2-a_1)\phi_1^2 - 2(b_1+b_2)\phi_1 -
      n\abs{d\phi_1}^2\right)\\ + \left( (a_1-a_2)\phi_2^2 -
      2(b_1+b_2)\phi_2- n\abs{d\phi_2}^2\right)
  \end{multline*}
  Since the left-hand side is constant and the first parenthesis on
  the right-hand depends only on the position $M_1$ whilst the second
  parenthesis only depends on the coordinate in $M_2$ the first
  statement follows.

  Obtaining the converse is a matter of inserting \eqref{eq:25} in
  \eqref{eq:21}.
\end{proof}

\subsection*{Proof of Theorem~\ref{thm:1}}

By formulas~\eqref{eq:6} and~\eqref{eq:17}, a product metric
$g=g_1+g_2$ on $M=M_1\times M_2$ is conformally Einstein if and only
if a positive function $\phi\colon M\to\mathbb R$ and a constant
$\bar\lambda = s(\bar g)/n$ exist such that
\begin{multline}
  \label{eq:7}
  \Ric(g_1) + \Ric(g_2) 
  =  - (n-2)\phi^{-1}\LC^g(d\phi) \\
  + (\bar\lambda \phi^{-2} + \phi^{-1}\Delta_g\phi +
  (n-1)\phi^{-2}\abs{d\phi}^2)(g_1+g_2).
\end{multline}
On the other hand, since $\bar g$ is Einstein it is mixed
Ricci-flat.  Therefore we get

\begin{lemma}\label{lem:3}
  If the Riemannian product $(M,g)$ is conformally Einstein and the
  summand $\phi_2$ in the conformal factor is non-constant then $g_1$
  is an Einstein metric and
  \begin{equation}\label{eq:29}
    \LC^{g_1}d\phi_1 = \frac1{n-2} \left(\bar\lambda\phi^{-1} -
      \lambda_1\phi + \Delta_g\phi +
      (n-1)\phi^{-1}\abs{d\phi}^2\right) g_1. 
  \end{equation}
\end{lemma}
\begin{proof}
  Insert $\phi = \phi_1 + \phi_2$ in equation~\eqref{eq:7} and project
  the result to $S^2(T^*M_1)$ to get
  \begin{equation} 
    \label{eq:18}
    \phi^2\Ric(g_1)=
    - (n-2)\phi\LC^{g_1}(d\phi_1) + (\bar\lambda + \phi\Delta_g\phi +
    (n-1)\abs{d\phi}^2)g_1.
  \end{equation}
  Differentiating this with respect to some $Y\in\Gamma(TM_2)$ gives
  \begin{multline}
    \label{eq:8}
    2\phi Y(\phi_2) \Ric(g_1) =  - (n-2)Y(\phi_2) \LC^{g_1}(d\phi_1)\\
    + \left(Y(\phi_2)\Delta_g\phi + Y(\Delta_2\phi_2)\phi +
      (n-1)Y(\abs{d\phi_2}^2\right)g_1.
  \end{multline}
  Suppose that $\phi_2$ is non-constant and pick $Y$ and $z$ in $M_2$
  such that $Y_z(\phi_2)\not= 0$.  Divide equation~\eqref{eq:8} by
  $Y(\phi_2)$, differentiate the result with respect to $Y$ and divide
  by $Y(\phi_2)$ again to get
  \begin{equation}
    \Ric(g_1) = \frac{Y(F)}{2Y(\phi_2)}g_1,\label{eq:24}
  \end{equation}
  where $F=\Delta_g\phi + \left(Y(\Delta_2\phi_2)\phi +
    (n-1)Y(\abs{d\phi_2}^2)\right)/Y(\phi_2)$.  Inserting
  (\ref{eq:24}) in (\ref{eq:18}) gives (\ref{eq:29}).
\end{proof}

If, say $\phi_2$, is constant set $B=M_1, ~F=M_2, ~g_B=g_2,
~f=\phi^{-1}$ and $g_B=f^2g_1$ to see that this is precisely the first
case of Theorem~\ref{thm:1}.  The conditions listed under \ref{item:1}
are both necessary and sufficient, see~\cite{Besse:Einstein}.  

So suppose now that both summands of the conformal factor are
non-constant.  By Lemma~\ref{lem:3}, $g_1$ and $g_2$ are both Einstein
and so $g$ has constant scalar curvature.  Taking the trace of
equation~\eqref{eq:29} and exploiting the symmetry of the setup leads
to
\begin{equation}\label{eq:31}
  -\frac{\Delta_1\phi_1}{n_1}+\frac{\lambda_1\phi}{n-2} = \frac1{n-2}
  \left(\bar\lambda\phi^{-1} + \Delta_g\phi +
    (n-1)\phi^{-1}\abs{d\phi}^2\right) =
  -\frac{\Delta_2\phi_2}{n_2} + \frac{\lambda_2\phi}{n-2}.
\end{equation}
Use equation~\eqref{eq:25} of Lemma~\ref{lem:2} in the equality
between far left and far right of~\eqref{eq:31} to get
\begin{gather*}
  \left(\frac{\lambda_1}{n-2} - \frac{a_1}{n_1}\right)\phi_1 +
  \frac{\lambda_1}{n-2}\phi_2 -\frac{b_1}{n_1} =
  \left(\frac{\lambda_2}{n-2} - \frac{a_2}{n_2}\right)\phi_2 +
  \frac{\lambda_2}{n-2}\phi_1 - \frac{b_2}{n_2}.
\end{gather*}
By the non-constancy of $\phi_1$ and $\phi_2$ this yields
\begin{gather}
  \label{eq:30}
  \frac{\lambda_1}{n-2} - \frac{a_1}{n_1} = \frac{\lambda_2}{n-2},\qquad
  \frac{\lambda_2}{n-2} - \frac{a_2}{n_2} =
  \frac{\lambda_1}{n-2},\quad\text{and}\quad 
  \frac{b_1}{n_1} = \frac{b_2}{n_2}.
\end{gather}
The equation $(n-1)(a_1+a_2)=s$ from Lemma~\ref{lem:2} combined
with~\eqref{eq:30} and some algebra ends with
\begin{equation}
  \label{eq:34}  
  s_1=(n_1-1)a_1,\qquad
  s_2=(n_2-1)a_2\quad\text{and}\quad n_2a_1+n_1a_2=0. 
\end{equation}
Setting $\bar a:= a_1/n_1$ and $\bar b:=b_1/n_1$ gives
part~\ref{item:2} and the claim about the Einstein constants of $g_1$
and $g_2$.

To see the converse of part~\ref{item:2} suppose $\bar a, ~\bar b$ to be
constants and $M_1$ and $M_2$ to be Einstein manifolds with
non-constant functions $\phi_1,~\phi_2$ satisfying
$\lc^{g_1}d\phi_1=-(\bar a\phi_1+\bar b)g_1$ and
$\lc^{g_2}d\phi_2=(\bar a\phi_1-\bar b)g_2$.  In this situation
$\lambda_1 = (n_1-1)\bar a$ and $\lambda_2 = -(n_2-1)\bar a$ since
this true near regular points for $\phi_j$ and since stationary points
are isolated, see~\cite{MR979791}.   We note that then $(n-2)\bar
a=\lambda_1-\lambda_2$.  On the other hand our hypothesis
leads to
\begin{align*}
  \phi\ric(g) + (n-2)\lc d\phi +\bar b g=& (\phi_1+\phi_2)(\lambda_1
  g_1+\lambda_2 g_2) + (n-2)((-\bar a\phi_1 g_1 + \bar a\phi_2 g_2)
  \\
  =& \left[(\lambda_1-(n-2)\bar a)\phi_1+\lambda_1\phi_2\right] g_1 \\&\quad+
  \left[\lambda_2\phi_1 + (\lambda_2+(n-2)\bar a)\phi_2\right]g_2,
\end{align*}
which is proportional to $g_1+g_2$ precisely when $\lambda_1- \lambda_2
= (n-2)\bar a$.  

Finally, the left hand equality in~\eqref{eq:31} is equivalent to
\begin{equation}\label{eq:32}
  0=(n-2)\phi\Delta_1\phi_1+n_1(\bar\lambda
  -\lambda_1\phi^2+\phi\Delta_g\phi + (n-1)\abs{d\phi}^2).
\end{equation}
Expand this according to equation~\eqref{eq:25},~\eqref{eq:30} and
\eqref{eq:34} and isolate the summands depending only on $M_1$ and
those depending only on $M_2$.  The sum of those depending on both is
zero.  Arguing as in the proof of Lemma~\ref{lem:2} we obtain
constants $\bar c_1$ and $\bar c_2$ such that $ \abs{d\phi_1}^2 + \bar
a \phi_1^2+ 2\bar b \phi_1 + \bar c_1= 0=\abs{d\phi_2}^2 - \bar a
\phi_2^2+ 2\bar b \phi_2 + \bar c_2$ and such that $\bar\rho=\bar c_1
+\bar c_2$.  This completes the proof of Theorem~\ref{thm:1}.\qed
\subsection*{Proof of Theorem~\ref{thm:2}}

Theorem~\ref{thm:2} follows by scrutinizing the results and
calculations collected in \cite{MR979791}.  We deal first with the
case $\bar a=0$.  

Given $M$ and $g$ as stated, Theorem~\ref{thm:1} tells us first that
$g$ is Ricci-flat and in particular Einstein.  Consulting the Main
Theorem of \cite{MR979791}, $(M,g)$ is Euclidean space and, checking
example $1b)$ of the same reference, $\phi(z)=\tfrac12(\abs{z}^2+R^2)$
for a positive constant $R$.  We conclude that $\bar b=-1$.  On the
other hand, $g_1$ is itself complete, Ricci-flat and admits a solution
to $\lc^{g_1}d\phi_1=g_1$ (and the same for $g_2$).  Invoke
Proposition~26 of \cite{MR979791} (due to Tashiro~\cite{MR0174022}) to
conclude that $M_1$ and $M_2$ are Euclidean spaces and $\phi_1(x) =
\tfrac12(\abs{x}^2+A_1)$, ~$\phi_2(x) = \tfrac12(\abs{x}^2+A_2)$ such
that $R^2 = A_1+A_2$.

Now suppose $\bar a=1$.  Theorem~25 of~\cite{MR979791} and its proof
tells us that $M_1$ is the standard sphere; writing $g_1 =
ds^2+\sin^2(s) g_{S^{n_1-1}}$ then, after some adjustment,
$\phi_1(s)=\cos(s)-\bar b$.  No equivalent to the Theorems 25 and 26
is available for solving $\lc^{g_2}d\phi_2=(\phi_2-\bar b)g_2$.
However, there are essentially only three possible solutions
distinguishable in the following way: if $\phi_2$ has a stationary
point then $M_2$ is hyperbolic space with metric expressed as $g_2 =
dt^2+\sinh^2(t)g_{S^{n_2-1}}$ and $\phi_2(t) = \cosh(t) + b$ since
$g_2$ is complete.  In this case $\phi=\cosh(t)+\cos(s)$, which has
zeros at the points $t=0$, $s=\pm\pi/2$.  For the other two cases
$M_2$ is a warped product globally $g_2=dt^2+(\phi_2')^2g_2^*$ where
$g_2^*$ is a complete Ricci-flat metric on some space $M_2^*$ and
$\phi_2$ is either $Ae^t+b$ or $A\cosh(t)+b$ for some constant $A$.
In both cases $\phi_2$ is globally defined on $M_2$ but
$\phi=\phi_1+\phi_2$ has zeros.

\qed

\begin{acknowledgements} The author was funded by the Junior Research
  Group ``Special Geometries in Mathematical Physics'' of the
  Volkswagen Foundation and the SFB 647, ``Space--Time--Matter'' of the
  DFG.
\end{acknowledgements}

\bibliographystyle{hamsplain}
\providecommand{\bysame}{\leavevmode\hbox to3em{\hrulefill}\thinspace}
\providecommand{\href}[2]{#2}

\end{document}